\documentclass[11pt]{article}
\usepackage{amsfonts}
\usepackage{subfigure} 
\usepackage{epsfig}
\usepackage{amsmath, amsthm}
\newcommand{\be}{\begin{eqnarray}}     	\newcommand{\ee}{\end{eqnarray}}

\newcommand{\vol}{\mathrm{Vol}}

\newcommand{\Int}{\mathrm{Int}}
\newcommand{\Hess}{\mathrm{Hess}}
 
\title{The rate of convergence of the mean curvature flow} 
\author{Tom Ilmanen, Natasa Sesum} 
\date{} 
\theoremstyle{plain}
\newtheorem{dummy}{Dummy} 

\theoremstyle{definition} 
\newtheorem{definition}[dummy]{Definition}
\newtheorem{step}{Step} 
\newtheorem{case}{Case} 
\theoremstyle{plain}
\newtheorem{corollary}[dummy]{Corollary}
\newtheorem{remark}[dummy]{Remark} 
\newtheorem{lemma}[dummy]{Lemma}
\newtheorem{theorem}[dummy]{Theorem}
\newtheorem{proposition}[dummy]{Proposition}

\newtheorem{claim}[dummy]{Claim} 

\begin{document}

\maketitle

\begin{abstract}
We study the flow $M_t$ of a smooth, strictly convex hypersurface by
its mean curvature in $\mathrm{R}^{n+1}$. The surface remains smooth
and convex, shrinking monotonically until it disappears at a critical
time $T$ and point $x^*$ (which is due to Huisken). This is equivalent
to saying that the corresponding rescaled mean curvature flow
converges to a sphere ${\bf S^n}$ of radius $\sqrt{n}$. In this paper
we will study the rate of exponential convergence of a rescaled
flow. We will present here a method that tells us the rate of the
exponential decay is at least $\frac{2}{n}$. We can define the
''arrival time'' $u$ of a smooth, strictly convex $n$-dimensional
hypersurface as it moves with normal velocity equal to its mean
curvature as $u(x) = t$, if $x\in M_t$ for $x\in \Int(M_0)$. Huisken
proved that for $n\ge 2$ $u(x)$ is $C^2$ near $x^*$. The case $n=1$
has been treated by Kohn and Serfaty, they proved $C^3$ regularity of
$u$. As a consequence of obtained rate of convergence of the mean
curvature flow we prove that $u$ is not $C^3$ near $x^*$ for $n\ge 2$.
We also show that the obtained rate of convergence $2/n$, that comes
out from linearizing a mean curvature flow is the optimal one, at
least for $n\ge 2$.
\end{abstract}

\begin{section}{Introduction}

In this paper we study a compact, smooth, strictly convex hypersurface
$M_0\in \mathrm{R}^{n+1}$ that moves with normal velocity equal to its
mean curvature. In other words, let $M_0$ be represented locally by a
diffeomorphism $F_0$ and let $F(\cdot,t)$ be a family of maps
satisfying the evolution equation
\begin{equation}
\label{equation-MCF_stan}
\frac{d}{dt}F = -H\nu,
\end{equation}
where $H(\cdot,t)$ is the mean curvature and $\nu(\cdot,t)$ is the
outer unit normal on $M_t$ and $M_t$ is the surface represented by
$F(\cdot,t)$. We often drop the $t$-dependence when no confusion will
result. Due to Huisken (see \cite{huisken1984}) the surface remains
smooth and convex and shrinks to a point. Assume it disappears at time
$T$ and that $x^*$ is a point to which it shrinks. Setting $x =
F(p,t)$, (\ref{equation-MCF_stan}) is then interpreted as
$$\frac{d}{dt}x = -H\nu (x).$$

The induced metric and the second fundamental form on $M$
will be denoted by $g = \{g_{ij}\}$ and $A = \{h_{ij}\}$. They can be
computed as follows:
$$g_{ij}(x) = \langle \frac{\partial F(x)}{\partial x_i},
\frac{\partial F(x)}{\partial x_j}\rangle,$$
$$h_{ij}(x) = -\langle \nu(x),\frac{\partial^2F(x)}{\partial x_i\partial x_j}
\rangle,$$
for $x\in \mathrm{R}^n$. The mean curvature is
$$H = g^{ij}h_{ij}.$$
We also use the notation
$$|A|^2 = g^{ij}g^{kl}h_{ik}h_{jl},$$
$$h = \frac{1}{\vol M}\int_M H dV.$$ 
In \cite{huisken1984} Huisken
computed the evolution equations of different curvatures.

\begin{theorem}[Corollary $3.5$ of \cite{huisken1984}]
\label{theorem-theorem_huisken1}
$$\frac{d}{dt}H = \Delta H + |A|^2H,$$
$$\frac{d}{dt}|A|^2 = \Delta |A|^2 - 2|\nabla A|^2 + 2|A|^4,$$
\begin{eqnarray*}
\frac{d}{dt}(|A|^2 - \frac{1}{n}H^2) &=& \Delta(|A|^2 - \frac{1}{n}H^2) 
- 2(|\nabla A|^2 - \frac{1}{n}|\nabla H|^2) \\
&+& 2|A|^2(|A|^2 - \frac{1}{n}H^2).
\end{eqnarray*}
\end{theorem}  

In order to prove his shrinking result (see Theorem $1.1$ in
\cite{huisken1984}) Huisken introduced a normalized flow, obtained by
reparametrization, keeping the total area of the evolving surface
fixed. He established important estimates for a normized flow
$$\frac{d}{d\bar{t}}\bar{F} = -\bar{H}\bar{\nu} +
\frac{1}{n}\bar{h}\bar{F},$$
where $\bar{F}(\cdot,t) = \psi(t)F(\cdot,t)$ and $\psi(t)$ is 
a function chosen so that the total area of $\bar{M}_t$ is being fixed and
$\bar{h} = \frac{1}{\bar{\vol}(\bar{M})}\int_{\bar{M}}\bar{H}^2$. Those
estimates are
\begin{equation}
\label{equation-equation_eq1'}
\bar{H}_{\max} - \bar{H}_{\min} \le Ce^{-\delta\bar{t}},
\end{equation}
\begin{equation}
\label{equation-equation_eq2'}
|\bar{A}|^2 - \frac{1}{n}\bar{H}^2 \le Ce^{-\delta\bar{t}},
\end{equation}
\begin{equation}
\label{equation-equation_eq3'}
|\nabla^m\bar{A}| \le C_me^{-\delta_m\bar{t}},\:\:\: m > 0,
\end{equation}
for some $\delta,\delta_m > 0$.  It is known that convex surfaces are
of type $1$ singularities (see \cite{huisken1984} and
\cite{huisken1990}). 

From now on, when we mention a rescaled flow, we will be thinking of
the following rescaling,
\begin{equation}
\label{equation-equation_our_scaling}
\tilde{F}(p,s) = (2(T-t))^{-1/2}F(p,t),
\end{equation}
with $s = -\frac{1}{2}\ln(T-t)$, where $T$ is a singularity time for
the original mean curvature flow. We will denote by
$\tilde{M}_{\tilde{t}}$ the rescaled surfaces moving by reparametrized
flow. The rescaled position vector then satisfies the equation
$$\frac{d}{dt}\tilde{F} = -\tilde{H}\nu + \tilde{F}.$$ 
In \cite{huisken1984} and \cite{huisken1990} Huisken showed that if the
expressions $P$ and $Q$, formed from $g$ and $A$, satisfy
$\frac{\partial P}{\partial t} = \Delta P + Q$ and if $\bar{P} =
\psi^{\alpha}P$ and if $\tilde{P} = (2(T-t))^{-\alpha/2}P$, then $\bar{Q}$
and $\tilde{Q}$ have degree $\alpha - 2$ and
$$\frac{d\bar{P}}{d\bar{t}} = \bar{\Delta}\bar{P} + \bar{Q} +
\frac{\alpha}{n}\bar{h}\bar{P},$$
$$\frac{d\tilde{P}}{d\tilde{t}} = \tilde{\Delta}\tilde{P} + \tilde{Q}
+ \alpha\tilde{P}.$$ 
Equations for $\bar{P}$ and $\tilde{P}$ look
quite similar and if one goes carefully through the estimates
established in \cite{huisken1984}, one can see that estimates
(\ref{equation-equation_eq1'}), (\ref{equation-equation_eq2'}),
(\ref{equation-equation_eq3'}) hold for corresponding quantities
$\tilde{A}$, $\tilde{H}$, etc. associated with rescaling
(\ref{equation-equation_our_scaling}). In particular this tells us
$\tilde{g}(s)$ uniformly converge to a round spherical metric, that
is, the surfaces $\tilde{M}_s$ are homothetic expansions of the
$M_t$'s and the surfaces $\tilde{M}_s$ converge to a sphere of radius
$\sqrt{n}$ in the $C^{\infty}$ topology as $s\to\infty$. 

\begin{remark}
The convergence of $\tilde{M}_s$ in any $C^k$-norm is exponential.
\end{remark}

We want to say more about this exponential convergence, that is, we
want to prove the following theorem.

\begin{theorem}
\label{theorem-rate-convergence}
If $M_0$ is uniformly convex, meaning that the eigenvalues of the
second fundamental form are strictly positive everywhere, then the
normalized equation (\ref{equation-tilda}) has a solution $x$ that
converges to a sphere of radius $\sqrt{n}$ exponentially at the rate
at least $\frac{2}{n}$.
\end{theorem}

Theorem \ref{theorem-rate-convergence} can be used to study the
arrival time of a smooth, strictly convex $n$-dimensional hypersurface
moving by a normal velocity equal to its mean curvature. Due to
Huisken we know the surface remains smooth and convex and shrinks to a
point $x^*$ at some finite time $T$. We define the ''arrival time'' on
the interior of the initial surface ($\partial\Omega = M_0)$ as $u(x)
= t$ if $x\in M_t$. A point $X^*$ to which a surface shrinks has a
unique maximum at $x^*$, $u(x^*) = T$. The smoothness of $u$ is
related to a roundness of $M_t$ as it shrinks to a point and it is the
best expressed in terms of the estimates for a curvature. Huisken
proved that $u$ is at least $C^2$ in $\Omega$ for $n\ge 2$.

The question whether $u$ is at least $C^3$ was raised by Kohn and
Serfaty in their recent work on a deterministic-control-based approach
to motion by curvature. Kohn and Serfaty proved that in the case
$n=1$, involving convex curves in the plane, $u$ is $C^3$ with
$D^3u(x^*) = 0$. The analogue of Huisken's work was done for curves in
the plane by Gage and Hamilton (see \cite{hamilton1986}). The
regularity of the arrival time was studied in this setting by Kohn and
Serfaty in \cite{kohn2004}. They needed at least $C^3$ regularity of
$u$ to draw a connection between a minimum exit time of two-person
game (see \cite{kohn2004} for more details) and the arrival time for a
curve shortening flow (see \cite{hamilton1986}). Their results would
completely extend to higher dimensions (drawing a connection between a
minimum exit time of the same game as above in higher dimensions and
the arrival time of a mean curvature flow) if we knew $u$ were $C^3$
near $x^*$. By Theorem \ref{theorem-rate-convergence} we can obtain
the following result.

\begin{theorem}
\label{theorem-regularity}
Function $u$ is not in general $C^3$ in $\Omega$ for $n\ge 2$.
\end{theorem}

We believe Theorem \ref{theorem-regularity} holds in the case $n=2$ as
well. In order to prove theorem \ref{theorem-regularity} we will
construct a solution to the rescaled mean curvature flow equation
whose behaviour is dictated by the first negative eigenvalue of an
operator $\Delta_{S^n} + 2$ (that is $-\frac{2}{n}$), which we obtain
while linearizing the rescaled mean curvature flow equation. As a
consequence of Theorem \ref{theorem-rate-convergence} and Theorem
\ref{theorem-regularity} we get the optimal rate of convergence 
of a mean curvature flow starting with a strictly convex
hypersurface.

\begin{theorem}
\label{theorem-rate-convergence}
The rate of convergence obtained in Theorem
\ref{theorem-rate-convergence} is the optimal one for $n\ge 2$.
\end{theorem}

The organization of the paper is as follows. In section 2 we derive a
linearization of a mean curvature flow equation. In section 3 we prove
the rate of exponential convergence of a strictly convex hypersurface
moving by (\ref{equation-MCF_stan}) is at least $2/n$, where $-2/n$
happens to be the biggest negative eigenvalue of a linear operator
$\Delta_{S^n} + 2$. In section 4 we prove Theorem
\ref{theorem-regularity} with a help of Theorem
\ref{theorem-rate-convergence}.  In section 5 we say more bout the
regualrity of $u$, that is, we give a condition on eigenvalues of
$\Delta_{S^n}+2$ (for $n\ge 2$) under which we can guarantee to
have some orders of regularity for $u$.

{\bf Acknowledgements:} I would like to thank R.Kohn and Sylvia
Serfaty for bringing the problem of regularity of $u$ to my attention
and for many useful discussions. 

\end{section}

\begin{section}{Linearizing the mean curvature flow equation around a sphere in 
$\mathrm{R}^{n+1}$}
\label{section-linearization}

In order to prove Theorem \ref{theorem-rate-convergence} we will study
a linearization of the rescaled mean curvature flow equation. It is a
standard matter, but for the sake of completeness we will include it
here.

\begin{definition}
A family of smoothly embedded hypersurfaces $(M_t)_{t\in I}$ in
$\mathrm{R}^{n+1}$ moves by {\it mean curvature} if
\begin{equation}
\label{equation-MCF}
\frac{d}{dt}x = -{\bf H}(x),
\end{equation}
for $x\in M_t$ and $t\in I$, $I\subset \mathrm{R}$ an open interval.
Here ${\bf H}(x)$ is the mean curvature vector at $x\in M_t$.
\end{definition}

Consider the family of smooth embeddings $F_t = F(\cdot,t): M^n\to
\mathrm{R}^{n+1}$, with $M_t = F_t(M^n)$ where $M^n$ is an
$n$-dimensional manifold. Setting $x=F(p,t)$, (\ref{equation-MCF}) is
then interpreted as
$$\frac{d}{dt}F(p,t) = -{\bf H}(F(p,t)),$$ 
for $p\in M^n$ and $t\in
I$. We can write a mean curvature vector as ${\bf H}(F(p,t)) =
H(p,t)\nu(p,t)$, where $H(\cdot,t)$ is the mean curvature and
$\nu(\cdot,t)$ is the outer unit normal on $M_t$.  We can define the
rescaled embeddings $\tilde{F}(p,s) = (2(T-t))^{-1/2}F(p,t)$, with
$s(t) = -\frac{1}{2}\ln(T-t)$. The surfaces 
$\tilde{M}_s = \tilde{F}(\cdot,s)(M)$ are defined for 
$-\frac{1}{2}\ln T \le s < \infty$ and satisfy the equation
$$\frac{d}{ds}\tilde{F}(p,s) = -\tilde{{\bf H}}(p,s) + \tilde{F}(p,s),$$
that is
\begin{equation}
\label{equation-tilda}
\frac{d}{ds}\tilde{x} = -\tilde{{\bf H}} + \tilde{x},
\end{equation}
if $\tilde{x} = \tilde{F}(p,s)$. In the rest of the paper we will be
considering evolution equation (\ref{equation-tilda}) and from now
on we will omit symbol $\:\: \tilde{}\:\:$ in symbols denoting the quantities 
characterizing the rescaled mean curvature flow. If we couple 
(\ref{equation-tilda}) with a normal $\nu$, we get
\begin{equation}
\label{equation-scalar}
\langle \frac{d}{dt}x, \nu\rangle = -H + \langle x,\nu\rangle.
\end{equation}
Consider the operator $L(x) = -H + \langle x,\nu\rangle$. We want 
to linearize it at a hypersurface given by $x$. In other words,
we want to compute $\frac{d}{ds}L(x_s)|_{s=0}$, where $x_s$ is a small
perturbation of $x$ (at some fixed time $t$) and $x_0 = x$. Let 
$u = \frac{d}{ds}x_s|_{s=0}$.

\begin{lemma}
\label{lemma-lin_normal}
$\frac{d}{ds}\nu(x_s)|_{s=0} = -\langle \nu,\frac{\partial u}{\partial
x_i}\rangle\frac{\partial F}{\partial x_j}g^{ij}$.
\end{lemma}

\begin{proof}
This is a straightforward computation: 
\begin{eqnarray*}
\frac{d}{ds}\nu(x_s)|_{s=0} &=& \langle \frac{d}{ds}\nu(x_s)|_{s=0},
\frac{\partial F}{\partial x_i}\rangle \frac{\partial F}{\partial x_j}g^{ij} \\
&=& -\langle \nu, \frac{\partial u}{\partial x_i}\rangle
\frac{\partial F}{\partial x_j}g^{ij}.
\end{eqnarray*}
\end{proof} 

\begin{lemma}
\label{lemma-lin_h}
$$\frac{d}{ds}h_{ij}|_{s=0} = \langle -\frac{\partial^2 u}{\partial
x_i\partial x_j} + \frac{\partial u}{\partial x_k}\Gamma_{ij}^k, \nu \rangle.$$
\end{lemma}

\begin{proof}
We know that the second fundamental form is given by a matrix
$$h_{ij} = -\langle \nu, \frac{\partial^2 F}{\partial x_i\partial
x_j}\rangle = \langle \frac{\partial F}{\partial
x_i},\frac{\partial\nu}{\partial x_j}\rangle.$$
We use Gauss-Weingarten relations
$$\frac{\partial^2 F}{\partial x_i\partial x_j} =
\Gamma_{ij}^k\frac{\partial F}{\partial x_k} - h_{ij}\nu,$$
to conclude 
\begin{eqnarray*} 
\frac{d}{ds}h_{ij}|_{s=0} &=& -\langle
\frac{\partial^2 u}{\partial x_i\partial x_j}, \nu\rangle +
\langle \frac{\partial^2 F}{\partial x_i\partial x_j}, 
\langle\nu,\frac{\partial u}{\partial x_k}\rangle
\frac{\partial F}{\partial x_l}g^{kl}\rangle \\
&=& -\langle\frac{\partial^2 u}{\partial x_i\partial x_j} + 
\frac{\partial u}{\partial x_k}\Gamma_{ij}^k, \nu\rangle,
\end{eqnarray*}
since $\langle \frac{\partial F}{\partial x_k},
\frac{\partial F}{\partial x_l}\rangle = g_{kl}$ and
$\langle \frac{\partial F}{\partial x_l}, \nu\rangle = 0$.
\end{proof}

\begin{lemma}
\label{lemma-lin_mean}
The linearization of the mean curvature $H$ is
$$-\frac{d}{ds}H|_{s=0} = \langle \Delta u, \nu\rangle + 
h_{ij}g^{ip}g^{jq}\{\langle \frac{\partial u}{\partial x_p}, 
\frac{\partial F}{\partial x_q}\rangle + \langle \frac{\partial u}{\partial x_q},
\frac{\partial F}{\partial x_p}\rangle\}.$$
where $u$ is a vector in $\mathrm{R}^{n+1}$ in a direction of a normal 
vector $\nu$.
\end{lemma}

\begin{proof}
Since $H = g^{ij}h_{ij}$, we have
\begin{equation}
\label{equation-medium_step}
\frac{d}{ds}H|_{s=0} = \frac{d}{ds}g_{ij}|_{s=0}h_{ij} + 
g^{ij}\frac{d}{ds}h_{ij}|_{s=0}.
\end{equation}
$$\frac{d}{ds}g_{ij}|_{s=0} = -g^{ip}g^{jq}\frac{d}{ds}g_{pq}|_{s=0},$$
and
$$\frac{d}{ds}g_{pq}|_{s=0} = \langle \frac{\partial u}{\partial x_p},
\frac{\partial F}{\partial x_q}\rangle + 
\langle \frac{\partial F}{\partial x_p}, \frac{\partial u}{\partial x_q}\rangle.$$
This together with (\ref{equation-medium_step}) and Lemma \ref{lemma-lin_h}
give
$$-\frac{d}{ds}H|_{s=0} = \langle \Delta u, \nu\rangle + 
h_{ij}g^{ip}g^{jq}\{\langle \frac{\partial u}{\partial x_p}, 
\frac{\partial F}{\partial x_q}\rangle + \langle \frac{\partial u}{\partial x_q},
\frac{\partial F}{\partial x_p}\rangle\}.$$
\end{proof}

\begin{proposition}
Let $u = \frac{d}{ds}x^s|_{s=0}$ where $x^s$ is a perturbation of $x$ and
let $w = \langle u, \nu \rangle$. Then
$$\frac{d}{dt}w = \Delta w + |A|^2w + w.$$
\end{proposition}

\begin{proof}
After taking $\frac{d}{ds}|_{s=0}$ of both sides of the evolution
equation $\langle\frac{d}{dt}x,\nu\rangle = -H + x\cdot\nu$ and using
Lemma \ref{lemma-lin_normal} and Lemma \ref{lemma-lin_mean} we get
\begin{eqnarray*}
\langle\frac{d}{dt}u,\nu\rangle + \langle \frac{d}{dt}x,-\langle
\nu,\frac{\partial u}{\partial x_i}\rangle\frac{\partial F}{\partial
x_j}g^{ij}\rangle &=& \langle \Delta u,\nu\rangle +
h_{ij}g^{ip}g^{jq}\{\langle \frac{\partial u}{\partial x_p},
\frac{\partial F}{\partial x_q}\rangle + \langle \frac{\partial u}{\partial x_q},
\frac{\partial F}{\partial x_p}\rangle\} +\\
&+& \langle u,\nu\rangle - \langle \nu,\frac{\partial u}{\partial x_i}\rangle g^{ij} 
x\cdot\frac{\partial F}{\partial x_j}.
\end{eqnarray*}
Since $\frac{d}{dt}x = -H\nu + x$, we have
\begin{equation}
\label{equation-scalar-w}
\langle\frac{d}{dt}u,\nu\rangle = \langle \Delta u,\nu\rangle +
h_{ij}g^{ip}g^{jq}\{\langle \frac{\partial u}{\partial x_p},
\frac{\partial F}{\partial x_q}\rangle + \langle \frac{\partial u}{\partial x_q},
\frac{\partial F}{\partial x_p}\rangle\} + \langle u,\nu\rangle.
\end{equation}

We will now compute $\Delta w$,
$$D_iw = \langle D_iu,\nu\rangle + \langle u, D_i\nu\rangle.$$
$$D_jD_i w = \langle D_jD_i u,\nu\rangle + \langle D_iu,D_j\nu\rangle
+ \langle D_ju, D_i\nu\rangle + \langle u,D_jD_i\nu\rangle.$$ 
By Gauss-Weingarten relation $\frac{\partial}{\partial x_j}\nu =
h_{jl}g^{lm}\frac{\partial F}{\partial x_m}$ we have
$$D_jD_i w = \langle D_jD_i u,\nu\rangle + h_{jl}g^{lm}\langle
D_iu,D_mF\rangle + h_{il}g^{lm}\langle D_ju, D_mF\rangle +
\langle u,D_jD_i\nu\rangle,$$
which gives
\begin{equation}
\label{equation-laplacian-w}
\Delta w = \langle \Delta u,\nu\rangle + \langle u,\Delta \nu\rangle
+ g^{ij}g^{lm}h_{jl}\langle D_iu,D_mF\rangle + 
g^{ij}g^{lm}h_{il}\langle D_ju,D_mF\rangle.
\end{equation}
Since $\frac{d}{dt}\nu = D^T H$, by (\ref{equation-scalar-w}) and 
(\ref{equation-laplacian-w}) we have
\begin{equation}
\label{equation-eq_w}
\frac{d}{dt}w = \langle u,D^T H\rangle + \Delta w - \langle
u,\Delta\nu\rangle + w.
\end{equation}

\begin{claim}
Let $M$ be a hypersurface in $\mathrm{R}^{n+1}$, given by an embedding
$F$. Then $\Delta\nu = D^TH - |A|^2\nu$, where $D^T H$ is a projection of
$DH$ onto a tangent space of a surface $M$.
\end{claim}

\begin{proof}
This is a pointwise computation, so we may assume $g_{ij}
=\delta_{ij}$ and $\Gamma_{ij}^k = 0$ at a particular point.  Denote
by $e_i = \frac{\partial f}{\partial x_i}$. Since $D_j\nu =
h_{jj}D_jF$, by Gauss-Weingarten relation $\frac{\partial^2
F}{\partial x_i\partial x_j} = \Gamma_{ij}^k\frac{\partial F}{\partial
x_k} - h_{ij}\nu$, we have $\langle D_j\nu, e_i\rangle = h_{ij}$,
and therefore
\begin{equation}
\label{equation-j}
\langle D_jD_j\nu,e_i\rangle + \langle D_j\nu,D_je_i\rangle =
D_jh_{ij} = D_ih_{jj},
\end{equation}
by Codazi equations. Moreover, $\langle D_j\nu,D_j e_i\rangle = 0$
at a point at which we are performing our computations. If we sum
all equations (\ref{equation-j}) over $j$ we get
\begin{equation}
\label{equation-ei}
\langle \Delta \nu, e_i\rangle = D^TH.
\end{equation}
By a similar computation we have
$$\langle D_jD_j\nu,\nu\rangle = -\langle D_j\nu,D_j\nu\rangle
= -h_{jl}h_{jm}g^{lm}.$$
If we sum the previous equations over $j$, we get
\begin{equation}
\label{equation-nu}
\langle \Delta \nu,\nu\rangle = -g^{jj}g^{lm}h_{lm}h_{jj} = -|A|^2.
\end{equation}
By our choice of coordinates at a point and by relations
(\ref{equation-ei}) and (\ref{equation-nu}) we have
$$\Delta\nu = D^TH - |A|^2\nu.$$
\end{proof}

The previous claim together with (\ref{equation-eq_w}) yield
\begin{equation}
\label{equation-scalar_version}
\frac{d}{dt}w = \Delta w + |A|^2w + w.
\end{equation}
\end{proof}

Let $x_{S^n}$ be an image of an embedding of a sphere ${\bf S^n}$ of
radius $\sqrt{n}$ into $\mathrm{R}^{n+1}$. Let $u = x - x_{S^n}$ and
$w = \langle u,\nu\rangle$.

\begin{lemma}
A scalar function $w$ satisfies the following evolution equation
$$\frac{d}{dt}w = \Delta_{S^n} w + 2w + Q,$$ 
where $\Delta_{S^n}$ is a Laplacian with respect to a metric on ${\bf S^n}$
and $Q$ is a quadratic term in $u$, $w$ and their first and second
covariant derivatives.
\end{lemma}

\begin{proof}
Since $x_{S^n}$ does not depend on time, it satisfies
$$\langle \frac{d}{dt}x_{S^n}, \nu\rangle = -H_{S^n} + \langle
x_{S^n},\nu_{S^n}\rangle,$$ 
because both sides of the previous
identity are equal to zero.  If we subtract this equation from
(\ref{equation-scalar}), we get
$$\langle \frac{d}{dt}(x-x_{S^n}), \nu\rangle = K(x) - K(x_{S^n}),$$
where $K(x) = -H(x) + x\cdot\nu$. By the previous consideration,
(\ref{equation-scalar_version}) and somewhat tedious, but standard computation,
we have
$$\frac{d}{dt}w = \Delta w + |A|^2w + w + Q,$$
where $Q$ is a quadratic term as in the statement of the lemma. 
Since $|A| = 1$ on a sphere of radius $\sqrt{n}$, we can
write the previous equation as
\begin{equation}
\label{equation-sphere}
\frac{d}{dt}w = \Delta_{S^n}w + 2w + Q',
\end{equation}
where $Q'$ is again a quadratic term in the same quantities as above, 
possibly different from $Q$.
Since $x(t)\to x_{S^n}$ exponentially, we can find a radial parametrization of
$M_t$ for $t$ sufficiently big, so that we can view $M_t$ as a radial graph 
over ${\bf S^n(\sqrt{n})}$ and consider $w$ as a scalar function 
defined on $S^n$.  
\end{proof}

\end{section}

\begin{section}{The rate of exponential convergence of the mean curvature flow}
\label{section-rate}

If $M_0$ is uniformly convex, i.e., the eigenvalues of its second
fundamental form are strictly positive everywhere. By results in
\cite{huisken1984} it follows that the rescaled equation
(\ref{equation-tilda}) has a solution that exponentially converges to
a sphere of radius $\sqrt{n}$. We want to say something more about the
rate of that exponential convergence. In order to do that we will
analyze the spectrum of $L(w) = \Delta_{S^n}w + 2w$. It is a standard
fact (see \cite{gallot}) that the spectrum of $L$ is given by
$\{-\frac{k(k+n-1)}{n} + 2\}_{k\in \{0\}\cup\mathrm{N}}$, if we adopt
the notation that $\Delta_{S^n}$ is a negative operator. The first
negative eigenvalue for $L$ is achieved for $k=2$ and is equal to
$-\frac{2}{n}$ (for $k = 0,1$ the corresponding eigenvalues are
$2,1$ respectively). This implies that $L$ does not have a zero
eigenvalue.

\begin{definition}
We will say that $x$ {\it converges} to a sphere $x_{S^n}$
exponentially {\it at a rate $\delta$} in $C^k$ norm, if there exist $C(k),t_0$
such that for all $t\ge t_0$,
$$|x-x_{S^n}|_k \le C(k)e^{-\delta t}.$$
\end{definition}

We may assume that $x$ converges to $x_{S^n}$ exponentially at rate
$\delta$. We want to say more about the rate of exponential
convergence.

\begin{theorem}
If $M_0$ is uniformly convex, meaning that the eigenvalues of the
second fundamental form are strictly positive everywhere, then the
normalized equation (\ref{equation-tilda}) has a solution $x$ that
converges to a sphere of radius $\sqrt{n}$ exponentially at the rate
at least $\frac{2}{n}$.
\end{theorem}

The ideas and techniques that we will use to prove Theorem
\ref{theorem-rate-convergence} rely on work of Cheeger and Tian (see
\cite{cheeger1994}). Similar approach has been used in \cite{natasa}
to prove the uniqueness of a limit of the Ricci flow.

Assume $\delta < 2/n$ since otherwise we are done. In order to prove
Theorem \ref{theorem-rate-convergence} we will use that the behaviour
of a solution of (\ref{equation-sphere}) is modeled on the behaviour
of a solution of a linear equation 
\begin{equation}
\label{equation-linear_eq}
\frac{d}{dt}v = L(v),
\end{equation}
where $L(v) = \Delta_{S^n}v + 2v$. Let $\{\lambda_k\}$ be the set of
all eigenvalues of $\mathcal{L}$. We can write $v = v_{\uparrow} +
v_{\downarrow} + v_0$, where $v_{\downarrow}(t) =
\sum_{\lambda_k<0}a_ke^{\lambda_kt}$, $v_{\uparrow}(t) =
\sum_{\lambda_k>0}a_k e^{\lambda_kt}$, and $v_0$ is a projection of
$v$ to a kernel of $L$. Since $L$ does not have zero
eigenvalue, $v_0 = 0$.

Some of the following consideration is based on the ideas and results
whose detailed proofs can be found in \cite{cheeger1994} (see also
\cite{natasa}) so we will just briefly outline them. The following
three lemmas can be found in \cite{cheeger1994} (see also
\cite{natasa}). The idea of considering three consecutive time
intervals is originally due to Simon (\cite{simon1983}).

\begin{lemma}
\label{lemma-lemma_alpha}
There exists $\alpha > 1$ such that
\begin{equation}
\label{equation-equation_growth}
\sup_{[K,2K]}||v_{\uparrow}|| \ge \alpha\sup_{[0,K]}||v_{\uparrow}||,
\end{equation}
\begin{equation}
\label{equation-equation_decay}
\sup_{[K,2K]}||v_{\downarrow}|| \le \alpha^{-1}\sup_{[0,K]}||v_{\downarrow}||.
\end{equation}
The norms considered above are standard $L^2$ norms. In particular, we can 
choose $\alpha = e^{\frac{2}{n}}$.
\end{lemma}

\begin{lemma}
\label{lemma-lemma_beta}
There exists $\beta < \alpha$ such that if
\begin{equation}
\label{equation-equation_impl1}
\sup_{[K,2K]}||v|| \ge \beta\sup_{[0,K]}||v||,
\end{equation}
then
\begin{equation}
\label{equation-equation_impl2}
\sup_{[2K,3K]}||v|| \ge \beta\sup_{[K,2K]}||v||,
\end{equation}
and if
\begin{equation}
\label{equation-equation_impl3}
\sup_{[2K,3K]}||v|| \le \beta^{-1}\sup_{[K,2K]}||v||,
\end{equation}
then
\begin{equation}
\label{equation-equation_impl4}
\sup_{[K,2K]}||v|| \le \beta^{-1}\sup_{[0,K]}||v||.
\end{equation}
Moreover, if $v_0 = 0$ at least one of
(\ref{equation-equation_impl2}), (\ref{equation-equation_impl4})
holds. If also $v_{\uparrow} = 0$, we can choose $\beta = e^{2/n}$.
\end{lemma}

The basic parabolic estimates (for example similarly as in
\cite{simon1983} and \cite{cheeger1994}) yield the following lemma.

\begin{lemma}
\label{lemma-lemma_regularity_est}
There exists $\tau > 0$ such that for any solution $w$ of
(\ref{equation-sphere}) with $|w(t_0)|_{k+2,\alpha} \le \tau$,
we have that
$$\sup_{(t_0,t_0+L)}|w(t)|_{k,\alpha} \le C\sup_{(t_0,t_0+L)}||w||,$$
where the first norm is $C^{k,\alpha}$ norm and the last norm is $L^2$
norm.
\end{lemma}

Let as in the previous section $u = x - x_{S^n}$ and $w = \langle u
,\nu\rangle$. It satisfies,
$$\frac{d}{dt}w = \Delta_{S^n}w + 2w + Q,$$ 
where $Q$ is a quadratic
term in $u$, $w$ and their first and second covariant derivatives. Let
$||\cdot||_{a,b} = \int_a^b|\cdot|$, where $|\cdot|$ is just the $L^2$
norm. Let $\pi$ denote an orthogonal projection on the subspace
$\ker(-\frac{d}{dt} + \Delta_{S^n} + 2)_{S^n(\sqrt{n})\times
[t_0,t_0+K]}$, with respect to norm $||\cdot||_{t_0,t_0+K}$. Put
$\pi w = (\pi w)_{\uparrow} + (\pi w)_{\downarrow}$.
The following proposition shows that the behaviour of
a solution of a linear equation (\ref{equation-linear_eq}) is modeled
on a behaviour of a solution of (\ref{equation-sphere}).
If $\epsilon > 0$ is any small number,  there is some $t_0$ so that
$|w(t)|_k < \epsilon$ for $t\ge t_0$.

\begin{proposition}
\label{proposition-proposition_comparison}
There exists $\epsilon_0 > 0$, depending on the uniform bounds on the
geometries $g(t)$, such that if $\epsilon < \epsilon_0$, then if
\begin{equation}
\label{equation-equation_sol1}
\sup_{[K,2K]}||w|| \ge \beta\sup_{[0,K]}||w||,
\end{equation}
then  
\begin{equation}
\label{equation-equation_sol2}
\sup_{[2K,3K]}||w|| \ge \beta\sup_{[K,2K]}||w||,
\end{equation}
and if
\begin{equation}
\label{equation-equation_sol3}
\sup_{[2K,3K]}||w|| \le \beta^{-1}\sup_{[K,2K]}||w||,
\end{equation}
then
\begin{equation}
\label{equation-equation_sol4}
\sup_{[K,2K]}||w|| \le \beta^{-1}\sup_{[0,K]}||w||,
\end{equation}
Moreover, since $(\pi w)_0 = 0$, at least one of
(\ref{equation-equation_sol2}), (\ref{equation-equation_sol4}) holds.
If $(\pi w)_{\uparrow} = 0$ we can choose $\beta = e^{2/n}$.
\end{proposition}

\begin{proof}
Assume there exist a sequence of constants $\tau_i\to 0$, and a
sequence of times $t_i\to\infty$ such that $|\eta_i(t)|_{k,\alpha} =
|w(t_i+t)|_{k,\alpha} \le \tau_i \to 0$ for all $t\ge 0$, but for
which none of the assertions in Proposition
\ref{proposition-proposition_comparison} holds. Let $\psi_i =
\frac{\eta_i}{\sup_{[K,2K]}|\eta_i|}$. Then in view of Lemma
\ref{lemma-lemma_regularity_est}, from standard compactness results
(see Lemma $5.22$ and Proposition $5.49$ in \cite{cheeger1994}) we get
that for a subsequence $\psi_i \stackrel{C^{k,\alpha'}}{\to} \psi$ and
$$\frac{d}{dt}\psi = \Delta_{S^n} \psi + 2\psi,$$
where $\psi$ has a property that contradicts Lemma
\ref{lemma-lemma_beta}. 
\end{proof}

\begin{lemma}
\label{lemma-no_growth}
If $v$ is a solution of (\ref{equation-linear_eq}), such that $|v| \le
Ce^{-\delta t}$, then $v_{\uparrow} = 0$.
\end{lemma}

\begin{proof}
If that is not the case, assume $v = \sum_{\lambda_k < 0}e^{\lambda_k
t} + be^{\gamma t} = \tilde{v} + be^{\gamma t}$, where $\gamma > 0$, $b
\neq 0$ and $\tilde{v} = \sum_{\lambda_k < 0}e^{\lambda_kt}$. By Lemma
\ref{lemma-lemma_alpha} we have that 
$$\sup_{[K,2K]}||\tilde{v}_{\downarrow}|| \le
\alpha^{-1}\sup_{[0,K]}||\tilde{v}_{\downarrow}||,$$
for $\alpha = e^{\frac{2}{n}}$. 
Applying Lemma \ref{lemma-lemma_beta} inductively to 
$\sup_{[iK,(i+1)K]}||\tilde{v}||$, for every $i$, we get
\begin{equation}
\label{equation-decay0}
||\tilde{v}|| \le Ce^{-\alpha t}.
\end{equation}
The fact that a rate of an exponentaial decay in (\ref{equation-decay0})
is given by the same $\alpha$ as in Lemma \ref{lemma-lemma_alpha} follows
immediatelly from the proof of Lemma $5.31$ in \cite{cheeger1994}.
Furthermore,
$$||be^{\gamma t}|| \le ||v|| + ||\tilde{v}|| \le
Ce^{-\min\{\delta,2/n\}t},$$
and we get a contradiction for big values of $t$ unless $b=0$. 
\end{proof}

We know that $w = \langle x-x_{S^n}, \nu\rangle$ solves the evolution
equation (\ref{equation-sphere}). Since $|w|_{k,\alpha} < Ce^{-\delta
t}$, by Lemma \ref{lemma-no_growth} we have $(\pi w)_{\uparrow} = 0$.
Since $(\pi w)_0 = 0$, by Proposition
\ref{proposition-proposition_comparison} at least one of
(\ref{equation-equation_sol2}), (\ref{equation-equation_sol4}) holds.
Since $w\to 0$ exponentially as $t\to\infty$, we have
(\ref{equation-equation_sol4}) holding with a rate of decay at least
$2/n$ because $(\pi w)_{\uparrow} = 0$. By using a parabolic regularity
theory we can get $C^k$ exponentail decay with the rate at least $2/n$.
We can now finish the proof of Theorem \ref{theorem-rate-convergence}.

\begin{proof}[Proof of Theorem \ref{theorem-rate-convergence}]
From the previous discussion we know that $|w|_k \le
C(k)e^{-\frac{2}{n}t}$. Let $t_0$ be such that $|w(t)|_{k+2} <
\epsilon$ for all $t\ge t_0$, where $\epsilon$ is taken from
Proposition \ref{proposition-proposition_comparison}. Assume that
$\gamma$ is a maximal rate of decay of $x$ to $x_{S^n}$, that is
$\gamma = \max\{\delta\:\:|\:\:$ exist $C$ such that $|w| \le
Ce^{-\delta t}, \:\:\forall t\ge t_0\}$. We may assume $\gamma < 2/n$,
since otherwise we are done.
\begin{equation}
\label{equation-perp}
\langle x-x_{S^n}, \nu_{S^n}\rangle = \langle
x-x_{S^n},\nu_{S^n}-\nu\rangle + \langle x-x_{S^n},\nu\rangle.
\end{equation}
Since $x\to x_{S^n}$ as $t\to\infty$ uniformly, we can regard $x$ as a
radial graph over $S^n$ and therefore $x-x_{S^n} \perp \nu_{S^n}$,
for $t$ sufficiently big, that is, $x-x_{S^n} = |x-x_{S^n}|\nu_{S^n}$.
From (\ref{equation-perp}) we get
$$|x-x_{S^n}| \le Ce^{-2\gamma t} + Ce^{-\frac{2}{n}t} =
Ce^{-\min\{2\gamma, 2/n\}t},$$
which contradicts the maximality of $\gamma$ unless $\gamma = 2/n$.
\end{proof}

\end{section} 

\begin{section}{Regularity of the arrival time function}

Due to Huisken (see \cite{huisken1993}) we know that the arrival time
function is at least of class $C^2$ in $\Omega = \Int(M_0)$. In the
case of $n=1$ (where instead of the mean curvature flow we deal with
the curve shortening flow) Kohn and Serfaty showed that $u$ is at
least $C^3$. The question that remains open is whether $u$ is $C^3$ or
more in higher dimensions ($n\ge 2$). It turns out it is not $C^3$
at $x^*$ in a generic case. 

Before we start proving Theorem \ref{theorem-regularity} lets first
slightly change the notation from above. Let $x$ satisfy
$$\frac{d}{dt}x = -\bar{H}\bar{\nu},$$
and $y = (2\tau)^{-1/2}x$, where $\tau = T - t$ and $s = -\frac{1}{2}\ln\tau$.
Quantities $\bar{H}$ and $\bar{\nu}$ correspond to the original mean curvature flow.
Then $y$ satisfies
$$\frac{d}{ds}y = -H\nu + y.$$
We may assume $y(s)$ converges to a sphere $S^n$ of radius $\sqrt{n}$.
We have derived in section \ref{section-linearization} that 
$w' = \langle y-y_{S^n},\nu\rangle$  satisfies
$$\frac{d}{ds}w' = \Delta_{S^n}w' + 2w' + Q(w'),$$
where $Q$ is a quadratic term in $w'$ and its first and second covariant derivatives.
Let $w = \langle y-y_{S^n},\nu_{S^n}\rangle$.

\begin{claim}
A scalar function $w$ satisfies
$$\frac{d}{ds}w = \Delta_{S^n}w + 2w + \tilde{Q}(w),$$
where $\tilde{Q}$ is an expression containing the quadratic terms in $w$, $w'$ and
their first and second covariant derivatives.
\end{claim}

\begin{proof}
\begin{eqnarray*}
\frac{d}{ds}w &=& \frac{d}{ds}w' + \frac{d}{ds}\langle y-y_{S^n},
\nu_{S^n}-\nu\rangle \\ 
&=& \Delta_{S^n}w' + 2w' + Qw' + \langle (-H\nu +
y) - (-H_{S^n}\nu_{S^n} + y_{S^n}), \nu_{S^n}-\nu\rangle + \\
&+&\langle y-y_{S^n}, -\nabla H + \nabla H_{S^n}\rangle \\
&=& \Delta_{S^n} w + 2w + \tilde{Q}(w,w').
\end{eqnarray*}
\end{proof}

Let $\mathcal{A}$ be a set of all solutions of
(\ref{equation-sphere}). Define a map $\psi: \mathcal{A}\to \mathrm{R}$
by $\psi({\bf a}) = \alpha$, where $\alpha$ is a coefficient of
$\phi e^{-\beta s}$ in $\pi{\bf a}$.

To prove Theorem \ref{theorem-regularity-again} we will need the
following Proposition that tells us how to construct solutions to the
rescaled mean curvature flow with a certain behaviour, dictated by the
first negative eigenvalue of $\Delta_{S^n} + 2$.

\begin{proposition}
\label{proposition-bad_solution}
There exists a solution $y$ to a rescaled mean curvature flow
(\ref{equation-tilda}) such that $\psi(\langle y - y_{S^n},
\nu_{S^n}\rangle) \neq 0$.
\end{proposition}

\begin{proof}
Fix a sphere $S^n$ of radius $\sqrt{n}$ and let $\phi$ be a
homogeneous, harmonic polynomial on $S^n$ corresponding to an
eigenvalue $-2/n$ of a differential operator $\Delta_{S^n} + 2$.
Consider a set of solutions $y_{\alpha}$ of
\begin{eqnarray*}
\frac{d}{ds}y_{\alpha} &=& -H\nu + y_{\alpha}, \\
y_{\alpha}(0) &=& \alpha\phi \nu_{S^n},
\end{eqnarray*}
for small values of $\alpha$, so that $y_{\alpha}$ is a strictly
convex hypersurface which by Huisken's result implies that every such
solution $y_{\alpha}(s)$ exponentially converges to a sphere of radius
$\sqrt{n}$. Let $w_{\alpha} = \langle y_{\alpha} - y_{S^n},
\nu_{S^n}\rangle$. We have seen that $w_{\alpha}$ satisfies
\begin{eqnarray}
\label{equation-help1}
\frac{d}{ds}w_{\alpha} &=& \Delta_{S^n}w_{\alpha} + 2w_{\alpha} +
Q(w_{\alpha}) \\
w_{\alpha}(0) &=& \alpha\phi. \nonumber
\end{eqnarray}
Our goal is to show there exists some $\alpha \neq 0$ so that a
solution $y_{\alpha}$ satisfies the property stated in the
proposition. The proof of the existence of such an $\alpha$ is given
in the few following steps.

\begin{step}
\label{step-comparable-small}
If $w$ is a solution to a nonlinear equation (\ref{equation-help1})
such that $|w|_l \le C(l)$ for all $l$ and all $s\in [0,L]$, and if
$k\ge 0$, there exist a uniform constant $C = C(L,k)$ and $\epsilon =
\epsilon(k)$, so that if $\sup_{s\in [0,L]}|w|_k < \epsilon$, then
$\sup_{s\in [0,L]}|w|_{k+1} < C\epsilon$.
\end{step}

\begin{proof}
The assertion tells us that if a solution is small in $C^k$ norm, it
will stay comparably small in $C^{k+1}$ norm. Assume without loss of
generality that $k=0$ (consideration for bigger $k$ is analogous). Our
goal is to show that $W^{2,l}$ norms of $w$ stay comparably small and
then to use Sobolev embedding theorems to draw the conslusion of the
assertion. If we multiply (\ref{equation-help1}) by $w$ and integrate
it over $S^n$,
$$\frac{d}{ds}\int w^2 + 2\int|\nabla w|^2 \le 2w^2 + Qw\cdot w.$$
Moreover, since $|Qw| \le C(|\nabla^2w||w| + |\nabla w|^2)$,
$$|\int Qw\cdot w| \le C\int|w|^2 + C\epsilon\int|\nabla w|^2,$$
and therefore, for small enough $\epsilon$, after integrating in 
$s\in [0,L]$, we get,
\begin{equation}
\label{equation-first-step}
\sup_{s\in [0,L]}\int w^2 (s) + 2\int_0^L\int|\nabla w|^2 \le
C\epsilon L = C(L)\epsilon,
\end{equation}
where we will use the same symbol $C$ to denote
different uniform constants and $C(L)$ to denote different
uniform constants depending on $L$. Apply a covariant derivative
(with respect to an induced metric on a sphere $S^n$) to
(\ref{equation-help1}), multiply it by $\nabla w$ and integrate
over $S^n$. A simple calculation yields
\begin{equation}
\label{equation-derivative}
\frac{d}{ds}\int|\nabla w|^2 + \int|\nabla^2w|^2 \le C\int|\nabla
w|^2 + \int \nabla(Q(w))*\nabla w,
\end{equation}
where we denote by $A*B$ any quantity obtained from $A\otimes B$
by one or more of the following operations: summation over 
pairs of matching upper and lower indices; contraction on
upper indices with respect to the metric; contraction on lower
indices with respect to the metric inverse; multiplication by 
uniform constants (\cite{dan}) or by uniformly bounded scalar 
functions (e.g. geometric quantities defined for $S^n$). Since
\begin{eqnarray*}
\int\nabla(Q(w))*\nabla w &=& - \int Q(w)*\nabla^2 w \\
&\le& C\epsilon \int|\nabla^2 w|^2 + C\int|\nabla w|^2,
\end{eqnarray*}
integrating (\ref{equation-derivative}) in $s$ and choosing $\epsilon$
small enough so that we can absorb $C\epsilon \int|\nabla^2 w|^2$ in
a corresponding term appearing on the right hand side of 
(\ref{equation-derivative}), we get
$$\sup_{s\in [0,L]}\int|\nabla w|^2 + \int_0^L\int|\nabla^2 w|^2
\le C\int_0^L\int|\nabla w|^2 \le C(L)\epsilon,$$
where we have used (\ref{equation-first-step}).
By taking more and more derivatives of (\ref{equation-help1}),
a similar consideration as above yields that
$$\sup_{s\in [0,L]}\int|\nabla^l w|^2(s) \le C(L,l,n)\epsilon.$$
By Sobolev embedding theorems we have that $|w|_1 \le C(L,n)\epsilon$,
and more general, if $\sup_{s\in [0,L]}|w|_k < \epsilon$, then
$\sup_{s\in [0,L]}|w|_{k+l} < C(L,l,n)\epsilon$.
\end{proof}

\begin{step}
\label{step-extension}
Fix $L > 0$. There exist $\epsilon = \epsilon(L,n)$ and $\delta =
\delta(L,n)$ so that if $|\alpha| < \epsilon$, then a solution
$w(s)$ exists for all $s\in [0,L]$ and $|w|_{C^0} < \delta$. 
\end{step}

\begin{proof}
A semigroup representation formula for $w_{\alpha}$ gives
$$w_{\alpha}(s) = e^{As}w_{\alpha}(0) + \int_0^s
e^{A(s-t)}Q(w_{\alpha}(t))dt,$$ 
where $A = \Delta_{S^n} + 2$. We will omit the subscript $\alpha$. The
spectrum of $A$ is given by $\{-\frac{k(k+n-1)}{n} +
2\}_{k=0}^{\infty}$. Denote those eigenvalues by $\lambda_k$. Denote
by $\psi_0$, $\psi_1$ the harmonic, homogeneous polynomials
(eigenfunctions of $A$) corresponding to $k=0$ and $k=1$,
respectively, by $\phi$ and $\phi_i$ the ones corresponding to $k=2$
and $k=i \ge 3$, respectively.  For every $\epsilon$ choose a maximal
time $\eta$ so that $w(s)$ exists for $s\in [0,\eta]$ and $|w(s)|_0 <
\delta$ (we will see how we choose $\delta$ later). We want to show 
that for small $\epsilon$ we can take $\eta$ to be at least $L$. 
Assume that it is not the case, that is, $\eta < L$ no matter how 
small $\epsilon$ we take. By Step
\ref{step-comparable-small} we get that $|w(s)|_2 < C(L)\delta$, for a
constant $C(L)$ that linearly depends on $L$, as we can see from
the consideration in the previous step. We can write 
$Q(w(s)) = \alpha_0(s)\psi_0 + \alpha_1(s)\psi_1 + \alpha_3(s)\phi + 
\sum_{k\ge 3} \alpha_1(s)\phi_k$, where $|Q(w(s))| \le C_1(L)\delta^2$,
for all $s\in [0,\eta]$, where $C_1(L)$ is now a constant 
that depends on $L$ quadratically. Then,
$$w(s) = \alpha\phi e^{-\frac{2s}{n}} +
\int_0^s(\alpha_0(t)e^{2(s-t)}\psi_0 + \alpha_1(t)e^{s-t}\psi_1 +
\sum_{i\ge 3}\phi_i\alpha_i(t)e^{\lambda_i(s-t)})dt.$$
Notice that all $\lambda_i < 0$. We have that for all $s\in [0,\eta]$,
$$|w(s)|_0 \le C\epsilon e^{-2L/n} + \delta^2 C_1(L)(Le^{2L} + Le^L + L).$$
Choose $\delta \le \frac{1}{3(C_1(L)L(2^{2L} + e^L))}$ and let 
$C\epsilon e^{-2L/n} < \delta/3$. Then,
$$|w(s)|_0 < \frac{2\delta}{3} < \delta,$$
which implies that for sufficiently small initial data 
(sufficiently small $\epsilon$) we can extend $w(s)$ beyond $\eta$
so that $|w(s)|_0 < \delta$ continues holding. This contradicts the
maximality of $\eta$, that is, for sufficiently small $\epsilon$
we have the conclusion of the step. 
\end{proof}

(*)Fix some big $3L$ and choose $\epsilon$ and $\delta$ as in Step
\ref{step-extension}. Our next goal is to show that for sufficiently
small $\epsilon, \delta > 0$ we can actually extend our solution $w$
(as a scalar function on $S^n$) all the way up to infinity, so that
$|w(s)|_0 < 2\delta$. For each small $\epsilon, \delta$ that satisfy
(*), find $L'$, that is, a maximal time so that $w(s)$ can be extended
all the way to $L'$, with $|w(s)|_0 < 2\delta$ holding. Subdivide
interval $[0,L']$ into subintervals of length $L$. We want to show
that for some choice of $\epsilon$, $L' = \infty$. Assume therefore
$L' < \infty$, no matter which choice for $\epsilon$ we make. Let
$\pi$ be as before, an orthogonal projection onto a subspace
$\ker(-\frac{d}{ds} + \Delta_{S^n} + 2)|_{S^n\times [iL,(i+1)L]}$ and
$w = (\pi w)_{\uparrow} + (\pi w)_{\downarrow}$.  Similarly as in
\cite{cheeger1994} we have that a behaviour of a solution of
(\ref{equation-help1}) is modeled on a behaviour of a solution of a
linear equation $\frac{d}{dt}F = \Delta_{S^n}F + 2F$ (see Proposition
\ref{proposition-proposition_comparison} in section
\ref{section-rate}).

\begin{step}
\label{step-no-growing-modes}
For sufficiently small $\epsilon$, where $|\alpha| < \epsilon$ we can
extend a solution $w_{\alpha}$ (call it only $w$) all the way up to
infinity so that $|w| < 2\delta$. 
\end{step}

\begin{proof}
Let $N$ be a maximal integer so that $[(N-1)L,NL]\subset [0,L']$ and
let $\gamma = \min\{|\lambda_i|\:\:|\:\:\lambda_i$ is an eigenvalue of
$\Delta_{S^n} + 2\}$. Since $(\pi w)_0 = 0$ on $S^n\times [(N-2)L,
(N-1)L]$, by Proposition \ref{proposition-proposition_comparison} we
have the following two cases.

\begin{case}
$\sup_{[(N-1)L,NL]}||w|| \ge e^{L\gamma/2} \sup_{[(N-2)L,(N-1)L]}||w||$.
\end{case}

This together with standard parabolic regularity imply
$$\sup_{[(N-2)L,(N-1)L]}|w|_k \le C(k)e^{-L\gamma/2}\delta \le \delta
e^{-L\gamma/3},$$ 
for $L$ sufficiently big, that we fix at the
beginning (from the previous estimate we see that its ''bigness''
depends on uniform constants; it is independent from the choices for
$\epsilon$ and $\delta$). By the same proof as in Step
\ref{step-extension}, that is, by our choice of $\epsilon$ and
$\delta$, considering $w((N-1)L)$ as an initial value, we get that $w$
can be extended to $[(N-1)L,(N+2)L]$ so that $|w(s)|_0 < 2\delta$ for
all $s\in [0,(N+2)L]$. To justify that, notice the following two
things: (a) as in Step \ref{step-extension} we can see that if
$\delta' < \delta$ we can choose smaller $\epsilon' < \epsilon$ so
that when $|w(s_0)| < \epsilon'$ then $\sup_{s\in [s_0,s_0+2L]}|w(s)|
< \delta'$, where everything is independent of the initial time $s_0$;
(b) by standard parabolic regularity, as in Step
\ref{step-comparable-small} we can get that $\sup_{s\in
[s_0,s_0+L]}||w(s)|| < 2\delta$ implies $\sup_{s\in
[s_0,s_0+L]}|w(s)|_k \le C(k)\delta$. This contradicts the 
maximality of $L'$, since $L' < (N+1)L$.

\begin{case}
$\sup_{[(N-1)L,NL]}||w|| \le e^{-L\gamma/2} \sup_{[(N-2)L,(N-1)L]}||w||$.
\end{case}

In this case (we may assume that $\delta$ is chosen so that $2\delta <
\eta$), applying Proposition \ref{proposition-proposition_comparison}
inductively, we get 
$$\sup_{[(N-1)L,NL]}||w|| \le e^{-NL\gamma/2}2\delta.$$ 
We can now argue similarly as in the previous case, that is we can again extend
solution $w$ past time $L'$ so that $|w(s)| < 2\delta$.
\end{proof}

This actually tells us there is an $\epsilon$ such that whenever
$|\alpha| < \epsilon$, then $L' = \infty$, that is, for sufficiently
small initial data, a solution $w$ to (\ref{equation-help1}) exists
and $|w(s)| < 2\delta$, where $\delta$ is taken to be small. By
Proposition \ref{proposition-proposition_comparison}, $w(s)$ has
either a growing or a decaying type of behaviour. If it had a growing
type of behaviour on some interval $[kL,(k+1)L]$, applying Proposition
\ref{proposition-proposition_comparison} inductively, we would get
that
$$C\delta > \sup_{[NL,(N+1)L]}||w(s)|| \ge
e^{(N-k)L\gamma/2}\sup_{[kL,(k+1)L]}||w(s)||,$$ 
for all $N$, which yields
a contradiction when $N\to\infty$.  In particular, this means that by
using the implication (\ref{equation-equation_sol3}) $\Rightarrow$ 
(\ref{equation-equation_sol4}) inductively and standard parabolic 
estimates we have that
$$|w(s)|_k \le C(k)\delta e^{-s\gamma /2},$$
for a uniform contant $C(k)$.

\begin{step}
\label{step-convergence}
There exists $\epsilon$, so that for $|\alpha| < \epsilon$, a solution
$y_{\alpha}(s)$ of a mean curvature flow
\begin{eqnarray}
\label{equation-MCF-again}
\frac{d}{ds}y_{\alpha} &=& -H\nu + y_{\alpha}, \\
y_{\alpha}(0) &=& \alpha \phi \nonumber,
\end{eqnarray}
converges exponentially to $S^n$ (the one that we have started with).
\end{step}

\begin{proof}
Let $\epsilon_0$ be such that whenever $|\alpha| \le \epsilon_0$, then
$y_{\alpha}(0)$ is a strictly convex hypersurface. We know in that
case $y_{\alpha}(s)$ converges in $C^k$ norm, exponentially, to a
sphere $S^n_{\alpha}$ of radius $\sqrt{n}$, and a quantity $\sup_{s\in
[0,\infty)}|y_{\alpha}(s) - y_{S^n}|_k$ makes sense. Define a function
$G: [-\epsilon_0,\epsilon_0] \to [0,\infty)$ by $G(\alpha) =
\sup_{s\in[0,\infty)}|y_{\alpha}(s) - y_{S^n}|_k$. It is a continuous
function and therefore bounded on a compact set $[0,\epsilon_0]$. This
implies all solutions $y_{\alpha}(s)$, for $\alpha\in [0,\epsilon_0]$
lie in a $C^k$ ball of a finite radius, with a center at
$y_{S^n}$. The continuity of this map implies that for sufficiently
small, say $\epsilon_0$, all solutions $y_{\alpha}$, for $\alpha\in
[-\epsilon_0,\epsilon_0]$ lie in a $C^0$ ball centred at $y_{S^n}$, of
radius $1/2$. This implies that every limit sphere $S^n_{\alpha}$ of a
solution $y_{\alpha}$ has a nonempty intersection with $S^n$.  It is a
well known result that if two solutions of a mean curvature flow
become disjoint, they stay disjoint for the remaining time of their
existence. That is why our solutions $y_{\alpha}(s)$ never become
disjoint from $S^n$ (**).

Denote by $w = \langle y - y_{S^n}, \nu_{S^n}\rangle$ (we actually
mean $y_{\alpha}$, but we are omitting the subscripts). An initial
hypersurface $y(0)$ can be written as an entire graph over $S^n$, that
is, for a choice of a unit normal $\nu$ for $M$, we have $\langle \nu,
\nu_{S^n}\rangle > 0$ everywhere on $M$. Choose $\alpha$ small
($|\alpha| < \epsilon$), as in Step \ref{step-no-growing-modes}, so
that an equation
\begin{eqnarray}
\label{equation-same-equation}
\frac{d}{ds}\tilde{w}(s) &=& \Delta_{S^n}\tilde{w} + 2\tilde{w} +
Q(\tilde{w}), \\
\tilde{w}(0) &=& \alpha \phi \nonumber,
\end{eqnarray}
has a solution all the way to infinity and $|\tilde{w}(s)|_k <
C(k)\delta e^{-s\gamma/2}$, where $C(k)$ is a uniform constant. Let
$\epsilon$ and $\delta$ be very small and let $\eta < \infty$ be a
maximal time such that $\langle \nu(s), \nu_{S^n}\rangle > 0$ for
$s\in [0,\eta)$. We can regard $w(s)$ as a function over $S^n$ for
$s\in [0,\eta)$, therefore satisfying
(\ref{equation-same-equation}). This implies $\tilde{w}(s) = w(s)$ and
$|w(s)|_k \le C(k)\delta e^{-s\nu/2}$, for $s\in [0,\eta)$. From
(\ref{equation-MCF-again}) we get,
\begin{equation}
\label{equation-inter}
\frac{d}{ds}\langle y-y_{S^n},\nu_{S^n}\rangle = -H\langle
\nu,\nu_{S^n}\rangle + \langle y,\nu_{S^n}\rangle.
\end{equation}
We have
\begin{eqnarray}
\label{equation-est1}
\langle y,\nu_{S^n}\rangle &=& \langle y-y_{S^n}, 
\nu_{S^n}\rangle + \sqrt{n} \nonumber\\
&=& w(s) + \sqrt{n} >\sqrt{n} - C\delta e^{-s\gamma/2},
\end{eqnarray}
and
\begin{eqnarray}
\label{equation-est2}
|\frac{d}{ds}w|_0 &=& |\Delta_{S^n}w + 2w + Q(w)|_0 \le C|w|_2 \nonumber \\
&\le& \tilde{C}\delta e^{-s\gamma/2}.
\end{eqnarray}
By (\ref{equation-inter}) and (\ref{equation-est1}),
$$\frac{d}{ds}w = \frac{d}{ds}\langle y - y_{S^n}, \nu_{S^n}\rangle >
-H\langle \nu, \nu_{S^n}\rangle + \sqrt{n} - C\delta e^{-s\gamma/2}.$$
Combining the last estimate together with (\ref{equation-est2}) 
yields,
$$H\langle \nu(s), \nu_{S^n}\rangle > \sqrt{n} - C\delta e^{-s\gamma/2} -
\tilde{C}\delta e^{-s\gamma/2}.$$
Our constants in the previous estimate are uniform and therefore 
if we make $\delta$ small enough (which we can achieve by decreasing
$\epsilon$), since $\gamma > 0$, we get
$$H\langle \nu(s), \nu_{S^n}\rangle > n/2,$$
for all $s\in [0,\eta)$. Since $H$ is bounded from above, we get
$$\langle \nu(s), \nu_{S^n}\rangle > \beta > 0,$$ 
for all $s\in [0,\eta)$. This implies 
the property $\langle \nu(s), \nu_{S^n}\rangle > 0$ continues 
holding for our solution $y(s)$ past time $\eta$, which
contradicts the maximality of $\eta$, unless $\eta = \infty$. This
together with (**) imply we can consider $w(s)$ as a function over
$S^n$ for all $s\in [0,\infty)$, satisfying
(\ref{equation-same-equation}). By uniqueness of solutions, we have
$w(s) = \tilde{w}(s)$ and henceforth $|w(s)| < Ce^{-s\gamma/2}$, for
all $s$, that is, $y_{\alpha}(s)$ converges exponentially to a sphere
$y_{S^n}$ when $|\alpha|$ is small.
\end{proof}

We can now finish the proof of Proposition
\ref{proposition-bad_solution}. Once we have a conclusion of Step
\ref{step-convergence}, similarly as in Lemma \ref{lemma-no_growth} we
can prove there are no growing modes in $w$, that is we can write
$$w(s) = \alpha\phi e^{-\frac{2s}{n}} + \int_0^s \sum_{i\ge
3}\phi_i\alpha_i(t)e^{\lambda_i(s-t)}dt,$$ 
where the notation is the same as in Step \ref{step-extension}.  
Similarly as in the Claim
\ref{claim-rate-R} we can show $\int_0^s \sum_{i\ge
3}\phi_i\alpha_i(t)e^{\lambda_i(s-t)}dt$ will decay at least at a rate
of $e^{-4s/n}$, so we can not expect any cancellations and since
$\pi(w)(0) = \alpha\phi$, we have $\pi w(s) = \alpha\phi e^{-2s/n}$,
where $\alpha$ is small, but can be taken to be different from zero.
\end{proof}

We can now finish the proof of Theorem \ref{theorem-regularity}.

\begin{proof}[Proof of Theorem \ref{theorem-regularity}]
Take a solution $y$ found by Proposition
\ref{proposition-bad_solution}. Since $y-y_{S^n} = (2\tau)^{-1/2}x -
y_{S^n}$, where $\tau = T - t = e^{-2s}$, we have that $w =
(2\tau)^{-1/2}\langle x,\nu_{S^n}\rangle - \sqrt{n}$. We know that
$y\to y_{S^n}$ exponentially and because of this uniform convergence
we can consider $y$ as a radial graph over $S^n$ for all sufficiently
big values of $s$. That is why we can write $|y| = \langle
y,\nu_{S^n}\rangle$ and $|x| = \langle x, \nu_{S^n}\rangle$ and
$$\tau = \frac{|x|^2}{2(w+\sqrt{n})^2},$$
which yields
$$u(x^*) - u(x) = \frac{|x|^2}{2n}(1 - 2\frac{w}{\sqrt{n}} +
3\frac{w^2}{n} + o(w^2)).$$
Let $\pi$ be a projection onto $\ker(-\frac{d}{ds} + \Delta_{S^n} +2)$
with $\pi(w)(0) = w(0)$.
Then $w = \pi w + R$, where $R$ is not in 
$\ker(-\frac{d}{ds} + \Delta_{S^n} +2)$ and
$$\pi w = \alpha\phi e^{-\beta s} + \sum_k\alpha_k\phi_k e^{-\beta_k
s},$$ 
where $\alpha \neq 0$ (justified by Proposition
\ref{proposition-bad_solution}), $\beta = 2/n$, $-\beta_k$ are the
remaining negative eigenvalues of $\Delta_{S^n} + 2$ and $\phi,\phi_k$ are
harmonic, homogenous polynomials restricted to a sphere
$S^n(\sqrt{n})$, of degrees $2$ and $k\ge 3$, respectively, and
$\alpha, \alpha_k$ are some constants. Because of Lemma
\ref{lemma-no_growth}, in $\pi w$ there are no contributions fom
eigenfuctions corresponding to positive eigenvalues of $\Delta_{S^n} +
2$. Then
$$\pi w = \alpha\phi\tau^{\beta/2} + o(\tau^{3\beta/2}).$$
We may assume $x^*$ is the origin in $\mathrm{R}^{n+1}$.
Since $\nabla u (x^*) = 0$ and $\nabla_i\nabla_j u(x^*) = -\frac{1}{n}\delta_{ij}$
(see \cite{huisken1993}), we have
\begin{equation}
\label{equation-tau}
\tau = u(0) - u(x) = \frac{1}{n}|x|^2 + O(|x|^3),
\end{equation}
which yields
$$u(0) - u(x) = \frac{|x|^2}{2n}(1 -
\alpha\frac{\phi}{\sqrt{n}}\frac{|x|^{\beta}}{n^{\beta/2}} + O(|x|^{3\beta/2}) + R).$$
Since $(2\tau)^{1/2} = \frac{|x|}{|y|}$, by (\ref{equation-tau}) we get
$$\frac{1}{2^{\beta/2}}\frac{|x|^{\beta}}{|y|^{\beta}} =
\frac{|x|^{\beta}}{n^{\beta/2}} + O(|x|^{3\beta/2}),$$
and therefore
\begin{equation}
\label{equation-asymptotic}
u(0) - u(x) = \frac{|x|^2}{2n}( 1 - 
\alpha\frac{\phi}{\sqrt{n}}\frac{1}{2^{\beta/2}}\frac{|x|^{\beta}}{|y|^{\beta}}
+ O(|x|^{3\beta/2}) + R).
\end{equation}

\begin{claim}
\label{claim-rate-R}
There is $\gamma \ge 2\beta$ so that $R = O(|x|^{\gamma})$.
\end{claim}

\begin{proof}
A scalar function $R$ satisfies 
\begin{equation}
\label{equation-R}
\frac{d}{ds}R = \Delta_{S^n}R + 2R + Qw.
\end{equation}
If $R = O(|x|^{\gamma}) = O(|y|^{\gamma}e^{-\gamma s})$, then
$\frac{d}{ds}R = O(|y|^{\gamma}e^{-\gamma s})$ and $\Delta_{S^n}R =
O(e^{-\gamma s})$. On the other hand, $Qw$ is of order $O(e^{-2\beta
s})$. Since $Q(w) \neq 0$, we have that $\gamma \ge 2\beta$, because
otherwise (\ref{equation-R}) could not be satisfied for large values
of $s$.
\end{proof}

So far we have found a solution $y(s)$ to a rescaled mean curvature
flow, whose existence, together with the asymptotic behaviour of its
arrival time given by (\ref{equation-asymptotic}), for $\alpha \neq 0$
is provided by Proposition \ref{proposition-bad_solution}. For $n\ge
3$, since $\beta = 2/n < 1$, from (\ref{equation-asymptotic}) it
follows immediatelly that $u(x)$ can not be $C^3$ at the origin.
In the case $n=2$ we have
$$u(0) - u(x) = \frac{|x|^2}{4} - \alpha
\frac{\phi(x)}{8\sqrt{2}}\frac{|x|}{|y|} + O(|x|)^{3/2 + 2}.$$
Take any $x$ such that $\phi(x) \neq 0$ and choose a line $tx$,
for $t\in \mathrm{R}$. Then, since $\phi$ is a homogeneous
polynomial of degree two,
$$u(0) - u(tx) = t^2\frac{|x|^2}{4} - 
\alpha t^2|t|\frac{\phi(x)}{8\sqrt{2}}\frac{|x|}{|y|} + O(|tx|)^{7/2}.$$
If we treat the right hand side as a function of $t$, we can see
that it is not $C^3$ at $t = 0$. Henceforth, $u(x)$ can not be $C^3$
at the origin.
\end{proof} 

In section \ref{section-rate} we have proved that a solution to a
rescaled mean curvature flow (\ref{equation-tilda}), starting with a
strictly convex hypersurface $M_0$, converges exponentially to a
sphere of radius $\sqrt{n}$ at a rate of at least $2/n$, that is
$$|y(s) - y_{S^n}| \le Ce^{-\delta s},$$
for $\delta \ge 2/n$. We will conclude that in a generic case we
can not expect to have $\delta > 2/n$.

\begin{proof}[Proof of Theorem \ref{theorem-rate-convergence}]
The proof goes by contradiction. Assume there is $\delta > 2/n$ so
that $\delta$ is the optimal rate of convergence of
(\ref{equation-tilda}) for $M_0$ being a strictly convex hypersurface.
Take a solution $y(s)$ that we have constructed in Proposition
\ref{proposition-bad_solution}. From the proof of the Proposition we
know that for some $\alpha \neq 0$, we have
\begin{equation}
\label{equation-contradiction}
\langle y(s) - y_{S^n}, \nu_{S^n}\rangle = \alpha\phi
e^{-\frac{2s}{n}} + \int_0^s \sum_{i\ge
3}\phi_i\alpha_i(t)e^{\lambda_i(s-t)}dt,
\end{equation} 
where $\lambda_i \le -(1+6/n)$ and we know that $y(s)$ converges to $y_{S^n}$
exponentially. If $|\langle y(s) - y_{S^n}, \nu_{S^n}\rangle| \le Ce^{-\delta s}$,
then (\ref{equation-contradiction}) would yield a contradiction for big 
values of $s$.
\end{proof} 

\end{section}

\begin{section}{More on regularity of $u(x)$ for some solutions
to the mean curvature flow}

If it happens that we have a solution $y$ such that $\psi(\langle
y-y_{S^n}, \nu_{S^n}\rangle) = 0$, our ''arrival time'' $u(x)$ might
be $C^3$ in $\Omega$. Moreover, the order of regularity depends on the
first term of form $\alpha_k\phi_ke^{-\beta_k s}$, appearing in
$\pi(\langle y-y_{S^n}, \nu_{S^n}\rangle)$, which actually determines
the rate of exponential convergence of $y$ to $y_{S^n}$. We will below
discuss the case of $C^3$-regularity, but the consideration is
analogous in the case of $C^k$-regularity, for $k > 3$.

\begin{corollary}
\label{corollary-regularity}
Let $y$ be a solution to (\ref{equation-tilda}) such that
$\psi(\langle y-y_{S^n}, \nu_{S^n}\rangle) = 0$ and $\pi w =
\sum_{k\ge l}\alpha_k\phi_ke^{-\beta_k s}$, for $l\ge 3$.  If $\beta_l
> 3$, then $u\in C^{3}(\Omega)$. This holds for any $n\ge 2$.
\end{corollary}

If our solution $y$ satisfies a condition in the corollary, as in
section \ref{section-rate}, we can prove that
$$|y - y_{S^n}|_k \le C(k) e^{-\beta_l s},$$
where $\beta_l \ge 1 + \frac{6}{n}$, that is, $y$ converges to a
sphere $y_{S^n}$ exponentially, at the rate at least $1 + \frac{6}{n}$.  

Function $u(x)$ can be viewed as the unique viscosity solution to the
nonlinear partial differential equation
$$\Delta u - \langle D^2u\frac{\nabla u}{|\nabla u|},\frac{\nabla
u}{|\nabla u|}\rangle + 1 = 0,$$
in $\Omega$ and $u=0$ at $\partial\Omega$. This equation was first 
studied by Evans and Spruck in \cite{evans1991}. They showed its solution
has the property that each level set $u=t$ is the smooth image of
$\partial\Omega$ under motion by curvature for time $t$, for any 
$0 \le t < T$. That is why the smoothness of $u(x)$ is apparent
away from $x^*$. Denote by $Z = |Du|^{-2}D_iuD_juD_iD_ju$. Then we can write 
the above equation as
\begin{equation}
\label{equation-elliptic-u}
\Delta u = Z - 1.
\end{equation}
Due to Huisken we know $u\in C^2$ with 
$\Hess_{ij} u = -\frac{1}{n}\delta_{ij}$. To prove $u\in C^3(\Omega)$ 
(for our flow $y$ having the properties as in the statement of the corollary)
we need to estimate $DZ$. Take $p\in \Omega$ and let $p\in M_t$, for some time 
$t$. We want to estimate $D_{\nu}Z(p)$ and $D_{\tau}Z(p)$, where $D_{\nu}Z$ is
a derivative of $Z$ in normal direction to the level set $M_t$ and $D_{\tau}Z$
is a tangential derivative at point $p\in M_t$. Teh estimate for $D_{\nu}Z$
is reduced to obtaining the estimate for $H^{-1}\frac{d}{dt}Z$, leading in
particular to  a term like $H^{-4}\Delta\Delta H$. All our hypersurfaces 
$M_t$ are embedded in $\mathrm{R}^{n+1}$ and for every function $f$ on
$\Omega$ we have that $(\nabla_{\mathrm{R}^{n+1}}f)^T = \nabla_{M_t}f$ at
$x\in M_t$. We will use $\nabla$ for $\nabla_{M_t}$. We need to estimate
$D_{\tau}Z$ which is translated to obtaining the estimate for 
$H^{-3} \nabla\Delta H$. If $\nu$ is the unit normal to $M_t$ then the
derivative of any function $f$ in the normal direction to the level set
$M_t$ of $u$ is given by $D_{\nu}f = H^{-1}\frac{d}{dt}f$. We can write
\begin{eqnarray*}
Z &=& D_{\nu}D_{\nu} u = D_{\nu}H^{-1} \\
&=& H^{-1}\frac{d}{dt}H^{-1} = -\frac{1}{H^3}\frac{dH}{dt} \\
&=& -\frac{1}{H^3}(\Delta H + |A|^2H) \\
&=& -\frac{1}{n} - (\frac{\Delta H}{H^3} + \frac{|A|^2 - \frac{1}{n}H^2}{H^2}).
\end{eqnarray*}

\begin{subsection}{Estimates on $D_{\nu}Z$ and $D_{\tau}Z$}

If we have $\frac{d}{dt}F = -H\nu$, in \cite{huisken1984} it was
computed that
$$\frac{d}{dt}g_{ij} = -2Hh_{ij}.$$ 
It is easy to compute the evolution equation  for the Christoffel
symbols (see \cite{dan})
$$\frac{d}{dt}\Gamma_{ij}^k = -g^{kl}(\nabla_i(Hh_{jl}) + \nabla_j(H
h_{il}) - \nabla_l(Hh_{ij})).$$
We  can compute
$$D_{\nu}Z = -\frac{\frac{d}{dt}\Delta H}{H^4} -
\frac{\frac{d}{dt}(|A|^2H)}{H^4} + 3\frac{(\Delta H +
|A|^2H)^2}{H^5}.$$
Fix $x\in M$ and a corresponding time $t_x$ such that $u(x) = t_x$.
Choose normal coordinates around $x$ in metric $g(t_x)$, so that 
$\Gamma_{ij}^k(x,t_x) = 0$. Since 
$\Delta H = g^{ij}(\nabla_i\nabla_j H - \Gamma_{ij}^k\nabla_k H)$,
\begin{eqnarray*}
\frac{d}{dt}\Delta H &=& 2g^{ip}g^{jq}Hh_{pq}(\nabla_i\nabla_j H -
\Gamma_{ij}^k\nabla_kH) + \Delta\Delta H + \\
&+& g^{ij}g^{kl}\nabla_kH(\nabla_i(Hh_{jl}) + \nabla_j(Hh_{il}) 
- \nabla_l(Hh_{ij})).
\end{eqnarray*}
By the curvature evolution equations (see
\ref{theorem-theorem_huisken1}) we get,
\begin{eqnarray}
\label{equation-Z}
D_{\nu}Z &=& -\frac{\Delta^2 H}{H^4} -
\frac{2g^{ip}g^{jq}h_{pq}\nabla_i\nabla_j H}{H^3} - \\
&-& \frac{g^{ij}g^{kl}\nabla_kH(\nabla_i(Hh_{jl}) + \nabla_j(Hh_{il})
- \nabla_l(Hh_{ij}))}{H^4} - \nonumber \\
&-& \frac{(\Delta|A|^2 + 2|\nabla A|^2)}{H^3}
+ 5\frac{|A|^2\Delta H}{h^4} + 3\frac{(\Delta H)^2}{H^5} \nonumber.
\end{eqnarray}
We want to discuss the asymptotics of terms appearing on the right hand 
side of identity (\ref{equation-Z}). 
\begin{itemize}
\item
Terms $\frac{\Delta^2 H}{H^4}$,
$\frac{2g^{ip}g^{jq}h_{pq}\nabla_i\nabla_j H}{H^3}$ can be estimated
by a constant times $(T-t)^{(\beta_l - 1)/2} \to 0$ as $t\to T$.
\item
Since the eigenvalues of $A$ are strictly positive, $|A|^2 \le H^2$ and  
$$|\frac{g^{ij}g^{kl}\nabla_kH}{H^3}| \le C\frac{|\nabla A|^2}{H^2} \le
C(T-t)^{(\beta_l-1)/2}.$$
\end{itemize}
We can similarly estimate the rest of the terms appearing in
(\ref{equation-Z}). The conclusion is that $|D_{\nu}Z| \le
C(T-t)^{\delta}$, for some $\delta > 0$.

As we have mentioned above we will use symbol $\nabla$ to denote a
derivative with respect to the induced metric $g(t)$ on $M_t$. We can
compute
$$\nabla Z = -\frac{\nabla\Delta H}{H^3} - \frac{\nabla(|A|^2H)}{H^4} +
3\frac{\Delta H\nabla H}{H^4} + 3\frac{\nabla H|A|^2}{H^3}.$$
All terms appearing on the right hand side of the above identity
are easy to estimate, e.g.
$$|\frac{\nabla\Delta H}{H^3}| \le C(T-t)^{(\beta_l - 1)}.$$
The conclusion is that $|D_{\tau}Z| \le C(T-t)^{\delta}$, for 
some $\delta > 0$. We can now finish the proof of Corollary
\ref{corollary-regularity}.

\begin{proof}[Proof of Corollary \ref{corollary-regularity}]
From (\ref{equation-elliptic-u}) we get
\begin{equation}
\label{equation-elliptic-Du}
\Delta Du = DZ.
\end{equation}
We know that 
$$|DZ (x)| \le C(T - t)^{(\beta_l -1)/2} = C(u(x^*) - u(x)) \le
\tilde{C}|x - x^*|^{(\beta_l -1)/2},$$
since $u\in C^1(\Omega)$. If $\beta_l > 3$ we get that $DZ$ is differentiable
at $x^*$, which means everywhere ($u$ is smooth everywhere on $\Omega$
except at $x^*$). This tells us $DZ$ lies in some H\"older space
$C^{0,\alpha}$, for $\alpha \in (0,1)$, and by elliptic regularity
applied to (\ref{equation-elliptic-Du}) we get $u\in C^{3,\alpha}(\Omega)$.
\end{proof}

\end{subsection}

\end{section}

\end{document}